\documentclass[preprint,12pt]{elsarticle}

\usepackage{amssymb}

\newtheorem{thm}{Theorem}

\newdefinition{dfn}{Definition}
\newdefinition{exm}{Example}
\newproof{pf}{Proof}

\newtheorem{lm}{Lemma}
\newtheorem{remark}{Remark}

\journal{Elsevier}

\begin{document}

\begin{frontmatter}



\title{On Some New Generalized Difference Sequence Spaces of Non-Absolute Type}

\author[label1]{Serkan Demiriz\corref{cor1}}
\ead{serkandemiriz@gmail.com}\cortext[cor1]{Corresponding Author
(Tel: +90 356 252 16 16, Fax: +90 356 252 15 85)}
\author[label2]{Osman \"{O}zdemir}
\ead{osman.ozdemir@gop.edu.tr}
\author[label3]{Osman Duyar}
\ead{osman5duy@hotmail.com}

\address[label1]{Department of Mathematics,  Gaziosmanpa\c{s}a University,
60250 Tokat, Turkey}
\address[label2]{Department of Mathematics,  Gaziosmanpa\c{s}a University,
60250 Tokat, Turkey}
\address[label3]{Anatolian High School, 60200 Tokat, Turkey \\
 }

\begin{abstract}
In this study, we define a new triangle matrix
$\widehat{W}=\{w_{nk}^{\lambda}(r,s,t)\}$ which derived by using
multiplication of $\lambda=(\lambda_{nk})$ triangle matrix with
$B(r,s,t)$ triple band matrix. Also, we introduce the sequence
spaces
$c_{0}^{\lambda}(\widehat{B}),c^{\lambda}(\widehat{B}),\ell_{\infty}^{\lambda}(\widehat{B})$
and $\ell_{p}^{\lambda}(\widehat{B})$   by using matrix domain of
this matrix on the sequence spaces $c_{0},c,\ell_{\infty}$ and
$\ell_{p}$ of the matrix $\widehat{W}$, respectively. Moreover, we
show that  norm isomorphic to the spaces $c_{0},c,\ell_{\infty}$ and
$\ell_{p}$, respectively. Furthermore, we establish some inclusion
relations concerning with those spaces and determine
$\alpha-,\beta-\gamma-$ duals of those spaces and construct their
Schauder basis. Finally, we characterize the classes
$(\mu_{1}^{\lambda}(\widehat{B}):\mu_{2})$ of infinite matrices ,
where $\mu_{1}\in\{c,c_{0},\ell_{p}\}$ and
$\mu_{2}\in\{\ell_{\infty},c,c_{0},\ell_{p}\}$.

\end{abstract}


\end{frontmatter}


\noindent
\section{Preliminaries,background and notation}
By a sequence spaces , we understand a linear subspace of the space
$\omega$  of all complex sequences which contains $\phi$, the set of
all finitely non-zero sequences.We write $\ell_\infty,c,c_0$ and
$\ell_p$ for the classical sequence spaces of all bounded,
convergent, null and absolutely $p$-summable sequences,
respectively, where $1\leq p<\infty$. Also by $bs$ and $cs$, we
denote the spaces of all bounded and convergent series,
respectively. We assume throughout unless stated otherwise that
$p,q>1$ with $p^{-1}+q^{-1}=1$ and use the convention that any term
with negative subscript is equal to zero. We denote throughout that
the collection of all finite subsets of $\mathbb{N}$ by
$\mathcal{F}$.

Let $A=(a_{nk})$ be  an infinite matrix of complex numbers $a_{nk}$,
where $n,k\in\mathbb{N}$. Then, $A$ defines a matrix mapping from
$X$ to $Y$ and is denote by $A:X\rightarrow Y$ if for every sequence
$x=(x_k)\in X$ the sequence $Ax=\{(Ax)_n\}$, the $A$-transform of
$x$, is in $Y$ where
\begin{equation}\label{e1.1}
(Ax)_n=\sum_{k} a_{nk}x_k, \quad   (n\in \mathbb{N})
\end{equation}
By $(X:Y)$, denote the class of all matrices $A$ such that
$A:X\rightarrow Y$. Thus, $A\in (X:Y)$ if and only if the series on
the right side (\ref{e1.1}) converges each $n\in \mathbb{N}$ and
$x\in X$, and we have $Ax=\{(Ax)_n\}_{n\in \mathbb{N}}\in Y$ for all
$x\in X$. A sequence $x\in \omega$ is said to be $A$-summable to $l$
if $Ax$ converges to $l$, which is called the $A$-limit of $x$.


A matrix $A=(a_{nk})$ is called a triangle if $a_{nk}=0$ for $k>n$
and $a_{nn}\neq 0$ for all $n,k\in\mathbb{N}$. It is trivial that
$A(Bx)=(AB)x$ holds for the triangle matrices $A,B$ and a sequence
$x$. Further, a triangle matrix $U$ uniquely has an inverse
$U^{-1}=V$ which is also triangle matrix. Then, $x=U(Vx)=V(Ux)$
holds for all $x\in\omega$.


Let us give the definition of some triangle limitation matrices
which are needed in the text. Let $q=(q_k)$ be a sequence of
positive reals and write
$$
Q_n=\sum_{k=0}^n  q_k, \quad (n\in \mathbb{N}).
$$
Then the Ces\`{a}ro mean of order one, Riesz mean with respect to
the sequence $q=(q_k)$ and Euler mean of order $r$ with $0<r<1$ are
respectively defined by the matrices $C=(c_{nk})$, $R^q=(r_{nk}^q)$
and $E^r=(e_{nk}^r)$; where
$$
c_{nk}=\left\{\begin{array}{ll}
  \displaystyle \frac{1}{n+1}, & (0 \leq k \leq n),\\
  0, & (k>n), \\
\end{array}\right.
\quad r_{nk}^{q}=\left\{\begin{array}{ll}
  \displaystyle \frac{q_{k}}{Q_{n}}, & (0 \leq k \leq n),\\
  0, & (k>n), \\
\end{array}\right.
$$
and
$$ e_{nk}^r=\left\{\begin{array}{cc}
  \displaystyle \left(n \atop k\right)(1-r)^{n-k}r^{k}, & (0 \leq k \leq n) \\
  0, & (k>n)
\end{array}\right.
$$
for all $k,n\in \mathbb{N}$. We write $\mathcal{U}$ for the set of
all sequences $u=(u_{k})$ such that $u_k\neq 0$ for all $k\in
\mathbb{N}$. For $u\in \mathcal{U}$, let $1/u=(1/u_k)$. Let
$z,u,v\in \mathcal{U}$, and define the summation matrix
$S=(s_{nk})$, the difference matrix $\Delta=(\Delta_{nk}^{(1)})$,
the generalized weighted mean or factorable matrix
$G(u,v)=(g_{nk})$, $A_u^r=\{a_{nk}^r(u)\}$,
$\Delta^{(m)}=(\Delta_{nk}^{(m)})$ by
$$
s_{nk}=\left\{\begin{array}{ll}
  \displaystyle 1, & (0 \leq k \leq n),\\
  0, & (k>n), \\
\end{array}\right.
\qquad \quad \Delta_{nk}^{(1)}=\left\{\begin{array}{ll}
  \displaystyle (-1)^{n-k}, & (n-1 \leq k \leq n),\\
  0, & (0\leq k<n-1\  \textrm {or}\ k>n), \\
\end{array}\right.
$$

$$
g_{nk}=\left\{\begin{array}{ll}
  \displaystyle u_{n}v_{k}, & (0 \leq k \leq n),\\
  0, & (k>n), \\
\end{array}\right.
\qquad \quad a_{nk}^{r}(u)=\left\{\begin{array}{ll}
  \displaystyle \frac{1+r^{k}}{n+1}u_{k}, & (0 \leq k \leq n),\\
  0, & (k>n), \\
\end{array}\right.
$$
and
$$
\Delta_{nk}^{(m)}=\left\{\begin{array}{cc}
  \displaystyle (-1)^{n-k}\left(m \atop n-k\right), & (\max\{0,n-m\}\leq k\leq n) \\
  0, & (0\leq k<n-1 \ \textrm{or} \ k>n)
\end{array}\right.
$$
for all $k,n \in \mathbb{N}$, where $u_n$ depends only on $n$ and
$v_k$ only on $k$. Let $r$ and $s$ be non-zero real numbers, and
define the generalized difference matrix $B(r,s)=\{b_{nk}(r,s)\}$ by
$$
b_{nk}(r,s)=\left\{\begin{array}{ll}
  \displaystyle r, & (k=n),\\
  \displaystyle s, & (k=n-1),\\
  0, & (0\leq k<n-1 \ \textrm{or} \ k>n),\\
\end{array}\right.
$$
for all $k,n \in \mathbb{N}$. The $B(r,s)$-transform of a sequence
$x=(x_k)$ is
$$
B(r,s)_k(x)=rx_k+sx_{k-1} \qquad \textrm{for all} \ k\in \mathbb{N}.
$$
We note that the matrix $B(r,s)$ can be reduced to the difference
matrices $\Delta$ in case $r=1$ and $s=-1$.


For a sequence space $X$, the matrix domain $X_A$ of an infinite
matrix $A$ is defined by
\begin{equation}\label{e1.2}
X_A=\{x=(x_k)\in \omega: Ax\in X\},
\end{equation}
which is a sequence  space. If $A$ is triangle, then one can easily
observe that the sequence space $X_A$ and $X$ are linearly
isomorphic, i.e., $X_A\cong X$.

Although in  the most cases the new sequence space $X_A$ generated
in the limitation matrix $A$ from a sequence space $X$ is the
expansion or the contraction of the original space $X$, it may be
observed in some cases that those space overlap. Indeed, one can
easily see that the inclusion, $X_S\subset X$ strictly holds for
$X\in \{\ell_{\infty},c,c_{0}\}$. As this, one can deduce that the
inclusion $X\subset X_{\Delta^{(1)}}$ also strictly holds for $X\in
\{\ell_{\infty},c,c_{0},\ell_{p}\}$. However, if we define
$X=c_0\oplus\textrm{span}\{z\}$ with $z=((-1))^{k}$, i.e., $x\in X$
if and only if $x:=\beta+\alpha z$ for some $\beta \in c_{0}$ and
some $\alpha \in \mathbb{C}$, and consider the matrix $A$ with the
rows $A_{n}$ defined by $A_{n}=(-1)^{n}e^{(n)}$ for all $n\in
\mathbb{N}$, we have $Ae=z\in X$ but $Az=e\notin X$ which lead us to
the consequences that $z\in X\backslash X_{A}$ and $e\in
X_{A}\backslash X$, where $e=(1,1,1,...)$ and $e^{(n)}$ is a
sequence whose only non-zero term is a $1$ in $n$th place for each
$n\in \mathbb{N}$. That is to say that the sequence spaces $X_{A}$
and $X$ overlap but neither contains to other.


The approach constructing a new sequence space by means of the
matrix domain of a particular limitation method has recently been
employed by Wang \cite{csw}, Ng and Lee \cite{nglee}, Malkowsky
\cite{ema1}, Altay and Ba\c{s}ar \cite{bafb1}, Malkowsky and
Sava\c{s} \cite{emaes}, Ba\c{s}ar{\i}r \cite{mb}, Ayd{\i}n and
Ba\c{s}ar \cite{cafb1}, Ba\c{s}ar et al. \cite{fbbamm},
\c{S}eng\"{o}n\"{u}l and Ba\c{s}ar \cite{msfb}, Altay \cite{ba},
Polat and Ba\c{s}ar \cite{hpfb}, and Malkowsky et al. \cite{emmmss}.
$\Delta, \Delta^{2}$ and $\Delta^{m}$ are the transposes of the
matrices $\Delta^{(1)}, \Delta^{(2)}$ and $\Delta^{(m)}$,
respectively, and $c_{0}(u,p)$ are the spaces consisting of the
sequences $x=(x_{k})$ such that $ux=(u_{k}x_{k})$ in the spaces
$c_{0}(p)$ and $c(p)$ for $u\in \mathcal{U}$, respectively, studied
by Ba\c{s}ar{\i}r \cite{mb}. More recently, the generalized
difference matrix $B(r,s)=\{b_{nk}(r,s)\}$ has been used by
Kiri\c{s}\c{c}i and Ba\c{s}ar \cite{mkfb} for generalizing the
difference spaces $\ell_{\infty}(\Delta), c(\Delta), c_{0}(\Delta)$,
and $bv_{p}$. Finally, the new technique for deducing certain
topological properties, for example $AB-, KB-, AD-$ properties,
etc., and determining the $\beta-$ and $\gamma-$ duals of the domain
of a triangle matrix in a sequence space has been given by Altay and
Ba\c{s}ar \cite{bafb5}.

Let $r,s$ and $t$ be non-zero real numbers, and define the
generalized difference matrix
$\widehat{B}=B(r,s,t)=\{b_{nk}(r,s,t)\}$ by
\begin{equation}\label{e1.3}
b_{nk}(r,s,t)=\left\{\begin{array}{ll}
  \displaystyle r, & (k=n) \\
  \displaystyle s, & (k=n-1) \\
  \displaystyle t, & (k=n-2) \\
  0, & (0\leq k<n-1 \ \textrm{or}\ k>n)
\end{array}\right.
\end{equation}
for all $n,k\in \mathbb{N}$. The inverse of
$B(r,s,t)=\{b_{nk}(r,s,t)\}$, which is denote
$B^{-1}(r,s,t)=\{d_{nk}(r,s,t)\}$ is given by
\begin{equation}\label{e1.4}
d_{nk}(r,s,t)=\left\{\begin{array}{ll}
  \displaystyle \frac{1}{r}\sum_{v=0}^{n-k}\bigg(\frac{-s+\sqrt{s^{2}-4tr}}{2r}\bigg)^{n-k-v}
\bigg(\frac{-s-\sqrt{s^{2}-4tr}}{2r}\bigg)^{v}, & (0 \leq k \leq n),\\
  0, & (k>n), \\
\end{array}\right.
\end{equation}
We should record here that $B(r,s,0)=B(r,s)$,
$B(1,-2,1)=\Delta^{(2)}$ and $B(1,-1,0)=\Delta^{(1)}$. So, the
results related to the matrix domain of the triple band matrix
$B(r,s,t)$ are more general and more comprehensive than the
consequences on the matrix domain of $B(r,s),\Delta^{(2)}$ and
$\Delta^{(1)}$, and  include them. We assume throughout that
$\lambda=(\lambda_{k})_{k=0}^{\infty}$ is a strictly increasing
sequence of positive reals tending to $\infty$, that is
\begin{equation}\label{e1.3}
0<\lambda_{0}<\lambda_{1}<... \qquad \textrm{and} \
\lim_{k\rightarrow \infty} \lambda_{k}=\infty.
\end{equation}

The main purpose of the present paper is to introduce  the sequence
space $\mu^{\lambda} (\widehat{B})$ and to determine the
$\alpha-,\beta-$ and $\gamma-$ duals of the space, where $\mu$
denotes the any of the classical spaces $\ell_{\infty},c,c_{0}$ or
$\ell_{p}$, and $\widehat{B}$ is the triple band matrix $B(r,s,t)$
and the sequence $\lambda=(\lambda_{k})$ is defined in (\ref{e1.3})
. Furthermore, the Schauder bases for the spaces
$c_{0}^{\lambda}(\widehat{B}),c^{\lambda}(\widehat{B})$ and
$\ell_{p}^{\lambda}(\widehat{B})$ are given, and some topological
properties of the spaces
$c_{0}^{\lambda}(\widehat{B}),c^{\lambda}(\widehat{B})$ and
$\ell_{p}^{\lambda}(\widehat{B})$ are examined. Finally, some
classes of matrix mappings  on the space $\mu^{\lambda}
(\widehat{B})$ are characterized.

The paper is organized as follows: In Section 2, the $BK-$ spaces
$c_{0}^{\lambda}(\widehat{B}),$ $c^{\lambda}(\widehat{B}),$
$\ell_{\infty}^{\lambda}(\widehat{B})$ and
$\ell_{p}^{\lambda}(\widehat{B})$ of generalized difference
sequences are introduced and  the Schauder bases of the spaces
$c_{0}^{\lambda}(\widehat{B}),c^{\lambda}(\widehat{B})$ and
$\ell_{p}^{\lambda}(\widehat{B})$ are given . In Section 3, are
examined some inclusion relations concerning with those space. In
Section 4, the $\alpha-,\beta-$ and $\gamma-$ duals of the
generalized difference sequence space $\mu^{\lambda} (\widehat{B})$
of non-absolute type is determined, respectively. In Section 5, the
classes $(\mu_{1}^{\lambda}(\widehat{B}):\mu_2)$ of infinite
matrices are characterized, where $\mu_{1}\in \{c,c_{0},\ell_{p}\}$
and $\mu_{2}\in \{\ell_{\infty},c,c_{0}, \ell_{p}\}$.
\section{The Difference Sequence Spaces $c_{0}^{\lambda}(\widehat{B}),c^{\lambda}(\widehat{B}),\ell_{\infty}^{\lambda}(\widehat{B})$
and $\ell_{p}^{\lambda}(\widehat{B})$ of Non-Absolute Type}

The difference sequence spaces have been studied by several authors
in different ways \cite{hk,cafb3, cafb4, cafb5, fbba, emamm,
zuamm,mas,carc,rcmeema,met,asfb}. In the present section, we
introduce the spaces
$c_{0}^{\lambda}(\widehat{B}),c^{\lambda}(\widehat{B}),\ell_{\infty}^{\lambda}(\widehat{B})$
and $\ell_{p}^{\lambda}(\widehat{B})$ and show that these spaces are
$BK-$ spaces of non-absolute type which are norm isomorphic to the
spaces $c_{0},c,\ell_{\infty}$ and $\ell_{p}$, respectively.
Furthermore, we give the basis for the spaces
$c_{0}^{\lambda}(\widehat{B}),c^{\lambda}(\widehat{B})$ and
$\ell_{p}^{\lambda}(\widehat{B})$.

Let $\lambda=(\lambda_{k})_{k=0}^{\infty}$ be a strictly increasing
sequence of positive reals tending to infinity, that is
$$
0<\lambda_{0}<\lambda_{1}<\lambda_{2}<... \quad  {\rm and} \quad
\lambda_{k}\rightarrow \infty \quad {\rm as} \quad k\rightarrow
\infty.
$$
We say that a sequence $x=(x_{k})\in \omega$ is $\lambda-$
convergent to the number $l\in \mathbb{C}$, called the $\lambda-$
limit of $x$, if $\Lambda_{n}(x)\rightarrow l$ as $n\rightarrow
\infty$ where
\begin{equation}\label{1.4}
\Lambda_{n}(x)=\frac{1}{\lambda_{n}}\sum_{k=0}^{n}
(\lambda_{k}-\lambda_{k-1})x_{k}; \quad (n\in \mathbb{N}).
\end{equation}
In particular, we say that $x$ is a $\lambda-$ null sequence if
$\Lambda_{n}(x)\rightarrow 0$ as $n\rightarrow \infty$. Further, we
say that $x$ is $\lambda-$ bounded if $\sup_{n\in \mathbb{N}}
|\Lambda_{n}(x)|<\infty$, \cite{mmakn1}. Recently, Mursaleen and
Noman \cite{mmakn1, mmakn2} studied the sequence spaces
$c_{0}^{\lambda},c^{\lambda}, \ell_{\infty}^{\lambda}$ and
$\ell_{p}^{\lambda}$ of non-absolute type as follows:
$$
c_{0}^{\lambda}=\bigg\{x=(x_{k})\in \omega: \lim_{n\rightarrow
\infty} \frac{1}{\lambda_{n}} \sum_{k=0}^{n}
(\lambda_{k}-\lambda_{k-1})x_{k}=0 \bigg\}
$$
$$
c^{\lambda}=\bigg\{x=(x_{k})\in \omega: \lim_{n\rightarrow \infty}
\frac{1}{\lambda_{n}} \sum_{k=0}^{n}
(\lambda_{k}-\lambda_{k-1})x_{k} \ \textrm{exists} \bigg\}
$$
$$
\ell_{\infty}^{\lambda}=\bigg\{x=(x_{k})\in \omega: \sup_{n\in
\mathbb{N}}\bigg| \frac{1}{\lambda_{n}} \sum_{k=0}^{n}
(\lambda_{k}-\lambda_{k-1})x_{k}\bigg|<\infty\bigg\}
$$
and
$$
\ell_{p}^{\lambda}=\bigg\{x=(x_{k})\in \omega:
\sum_{n=0}^{\infty}\bigg|\frac{1}{\lambda_{n}} \sum_{k=0}^{n}
(\lambda_{k}-\lambda_{k-1})x_{k}\bigg|^{p}<\infty\bigg\}.
$$
On the other hand, we define the matrix
$\widehat{\Lambda}=(\widehat{\lambda}_{nk})$ for all $n,k\in
\mathbb{N}$ by
\begin{equation}\label{e2.2}
\widehat{\lambda}_{nk}=\left\{\begin{array}{ll}
  \displaystyle
  \frac{\lambda_{k}-\lambda_{k-1}}{\lambda_{n}}, & (0 \leq k \leq n) \\
  0, & (k>n).
\end{array}\right.
\end{equation}
Then, it can be easily seen that the equality
\begin{equation}
\widehat{\Lambda}_{n}(x)=\frac{1}{\lambda_{n}}\sum_{k=0}^{n}
(\lambda_{k}-\lambda_{k-1})x_{k}
\end{equation}
holds for all $n\in \mathbb{N}$ and every $x=(x_{k})\in \omega$,
which leads us together with (\ref{e1.2}) to the fact that
$$
c_{0}^{\lambda}=\{c_{0}\}_{\widehat{\Lambda}}, \quad
c^{\lambda}=c_{\widehat{\Lambda}}, \quad
\ell_{\infty}^{\lambda}=\{\ell_{\infty}\}_{\widehat{\Lambda}}, \quad
\ell_{p}^{\lambda}=\{\ell_{p}\}_{\widehat{\Lambda}}.
$$

More recently, S\"{o}nmez \cite{asfb2} has defined the sequence
spaces $\ell_{\infty}(\widehat{B}), c(\widehat{B}),
c_{0}(\widehat{B})$ and $\ell_{p}(\widehat{B})$ as follows:
\begin{eqnarray*}
\ell_{\infty}(\widehat{B})&=&\bigg\{x=(x_{k})\in \omega: \sup_{k\in
\mathbb{N}}|rx_{k}+sx_{k-1}+tx_{k-2}|<\infty \bigg\},\\
c(\widehat{B})&=&\bigg\{x=(x_{k})\in \omega: \exists l\in \mathbb{C}
\ni \lim_{k\rightarrow \infty
 }|rx_{k}+sx_{k-1}+tx_{k-2}-l|=0  \bigg\},\\
c_{0}(\widehat{B})&=&\bigg\{x=(x_{k})\in \omega: \lim_{k\rightarrow
\infty
 }|rx_{k}+sx_{k-1}+tx_{k-2}|=0 \bigg\},\\
\ell_{p}(\widehat{B})&=&\bigg\{x=(x_{k})\in \omega: \sum_{k}
|rx_{k}+sx_{k-1}+tx_{k-2}|^{p}<\infty\bigg\}.
\end{eqnarray*}
In fact, the sequence spaces $\ell_{\infty}(\widehat{B}),
c(\widehat{B}), c_{0}(\widehat{B})$ and $\ell_{p}(\widehat{B})$ can
be consider as the set of all sequences whose $B(r,s,t)-$ transforms
are in the spaces $\ell_{\infty},c,c_{0}$ and $\ell_{p}$,
respectively. That is,
$$
\ell_{\infty}(\widehat{B})=\{\ell_{\infty}\}_{B(r,s,t)}, \quad
c(\widehat{B})=c_{B(r,s,t)}, \quad
c_{0}(\widehat{B})=\{c_{0}\}_{B(r,s,t)}, \quad
\ell_{p}(\widehat{B})=\{\ell_{p}\}_{B(r,s,t)}.
$$

Now, we introduce  the difference sequence spaces
       $\ell_{\infty}^{\lambda}(\widehat{B}),c^{\lambda}(\widehat{B}), c_{0}^{\lambda}(\widehat{B})$
and $\ell_{p}^{\lambda}(\widehat{B})$ as follows:
\begin{eqnarray*}
\ell_{\infty}^\lambda(\widehat{B})&=&\bigg\{x=(x_{n})\in
\omega:\sup_{n\in \mathbb{N}}
 |\frac{1}{\lambda_n}\sum_{k=0}^n(\lambda_k-\lambda_{k-1})(rx_{k}+sx_{k-1}+tx_{k-2})|<\infty \bigg\},\\
c^\lambda(\widehat{B})&=&\bigg\{x=(x_{n})\in \omega:
\lim_{n\rightarrow \infty
 }\frac{1}{\lambda_n}\sum_{k=0}^n(\lambda_k-\lambda_{k-1})(rx_{k}+sx_{k-1}+tx_{k-2})\ \textrm{exists} \bigg\},\\
c_{0}^\lambda(\widehat{B})&=&\bigg\{x=(x_{n})\in \omega:
\lim_{n\rightarrow \infty
 }\frac{1}{\lambda_n}\sum_{k=0}^n(\lambda_k-\lambda_{k-1})(rx_{k}+sx_{k-1}+tx_{k-2})=0 \bigg\},\\
\ell_{p}^\lambda(\widehat{B})&=&\bigg\{x=(x_{n})\in \omega: \sum_{n}
\bigg|\frac{1}{\lambda_n}\sum_{k=0}^n(\lambda_k-\lambda_{k-1})(rx_{k}+sx_{k-1}+tx_{k-2})\bigg|^{p}<\infty\bigg\}.
\end{eqnarray*}


On the other hand, we define the triangle matrix
$\widehat{W}=\{w_{nk}^{\lambda}\}=\widehat{\Lambda}\widehat{B}$ by
\begin{equation}\label{e2.3}
w_{nk}^{\lambda}=\left\{\begin{array}{ll}
 \displaystyle \frac{r(\lambda_{k}-\lambda_{k-1})+s(\lambda_{k+1}-\lambda_{k})+t(\lambda_{k+2}-\lambda_{k+1})}{\lambda_{n}}, & (k<n-1)\\
  \displaystyle \frac{r(\lambda_{n-1}-\lambda_{n-2})+s(\lambda_{n}-\lambda_{n-1})}{\lambda_{n}}, & (k=n-1) \\\
  \displaystyle \frac{r(\lambda_{n}-\lambda_{n-1})}{\lambda_{n}}, & (k=n) \\0, & (\textrm {otherwise})
\end{array}\right.
\end{equation}
for all $k,n\in \mathbb{N}$. Then , it can be easily seen that the
equality
\begin{equation}\label{e2.7}
\widehat{W}_{n}(x)=\frac{1}{\lambda_{n}}\sum_{k=0}^{n}
(\lambda_{k}-\lambda_{k-1})(rx_{k}+sx_{k-1}+tx_{k-2})
\end{equation}
holds for all $n\in \mathbb{N}$ and every $x=(x_k)\in \omega$. In
fact, the sequence spaces
$c_{0}^{\lambda}(\widehat{B}),c^{\lambda}(\widehat{B}),\ell_{\infty}^{\lambda}(\widehat{B})$
and $\ell_{p}^{\lambda}(\widehat{B})$ can be consider as the set of
all sequences whose $\widehat{W}-$ transforms are in the spaces
$c_{0},c,\ell_{\infty}$ and $\ell_{p}$, respectively. That is,
\begin{equation}\label{e11}
c_{0}^{\lambda}(\widehat{B})=\{c_0\}_{\widehat{W}},\quad
c^{\lambda}(\widehat{B})=c_{\widehat{W}},\quad
\ell_\infty^{\lambda}(\widehat{B})=\{\ell_\infty\}_{\widehat{W}},\quad
\ell_p^{\lambda}(\widehat{B})=\{\ell_p\}_{\widehat{W}}.
\end{equation}


Further, for any sequence $x=(x_k)$ we define the sequence
$y_k(\lambda)=\{y_k(\lambda)\}$ which will be used, as the
$\widehat{W}$-transform of $x$ and so we have
\begin{eqnarray}\label{e1..1}
y_{k}(\lambda)&=&\sum_{j=0}^{k-2}\frac{r(\lambda_{j}-\lambda_{j-1})+s(\lambda_{j+1}-\lambda_{j})+t(\lambda_{j+2}-\lambda_{j+1})}{\lambda_k}x_j\\
&+&\frac{r(\lambda_{k-1}-\lambda_{k-2})+s(\lambda_{k}-\lambda_{k-1})}{\lambda_k}x_{k-1}+\frac{r(\lambda_{k}-\lambda_{k-1})}{\lambda_k}x_{k};
\quad (k\in \mathbb{N}) \nonumber.
\end{eqnarray}

Since the proof may also be obtained in the similar way as for the
other spaces, to avoid the repetition of the similar statements, we
give the proof only for one of those spaces. Now, we may begin with
the following theorem which is essential in the study.

\begin{thm}\label{t1}
(i) The difference sequence spaces $c_{0}^{\lambda}(\widehat{B})$,
$c^{\lambda}(\widehat{B})$ and $\ell_\infty^\lambda(\widehat{B})$
are $BK-$ space with the norm
$\|x\|_{c_{0}^{\lambda}(\widehat{B})}=\|x\|_{c^{\lambda}(\widehat{B})}=\|x\|_{\ell_\infty^\lambda(\widehat{B})}
=\|\widehat{W}(x)\|_\infty $, that is ,
$$
\|x\|_{c_{0}^{\lambda}(\widehat{B})}=\|x\|_{c^{\lambda}(\widehat{B})}=\|x\|_{\ell_\infty^\lambda(\widehat{B})}
=\sup_{n\in\mathbb{N}}\big|\widehat{W}_n(x)\big|.
$$

(ii) Let $1\leq p<\infty$. Then $\ell_{p}^{\lambda}(\widehat{B})$ is
a $BK-$ space with the norm
$\|x\|_{\ell_{p}^{\lambda}(\widehat{B})}=\|\widehat{W}x\|_{p}$, that
is,
$$
\|x\|_{\ell_{p}^{\lambda}(\widehat{B})}=\bigg(\sum_{n}
|\widehat{W}_{n}(x)|^{p}\bigg)^{1/p}.
$$
\end{thm}

\begin{pf}
Since (\ref{e11}) holds and $c_{0},c$ and $\ell_{\infty}$ are $BK-$
spaces with respect to their natural norms (see \cite[pp.
16-17]{fb}) and the matrix $\widehat{W}$ is a triangle, Theorem
4.3.12 Wilansky \cite[pp. 63]{aw} gives the fact that
$c_{0}^{\lambda}(\widehat{B})$, $c^{\lambda}(\widehat{B})$ and
$\ell_\infty^\lambda(\widehat{B})$ are $BK-$ spaces with the given
norms. This completes the proof .
\end{pf}

\begin{remark}
Let $\mu\in\{\ell_{\infty},c,c_{0},\ell_{p}\}$ and  $|x|=|x_k|$.
Then the absolute property does not hold on the space
$\mu^{\lambda}(\widehat{B})$, that is,
$\|x\|_{\mu^{\lambda}(\widehat{B})}\neq
\||x|\|_{\mu^{\lambda}(\widehat{B})} $. This can be shown that for
at least one sequence in those space. Hence
$\mu^{\lambda}(\widehat{B})$ is the sequence space of non-absolutely
type.
\end{remark}

\begin{thm}\label{t2}
The sequence spaces $c_{0}^{\lambda}(\widehat{B})$,
$c^{\lambda}(\widehat{B})$, $\ell_{\infty}^{\lambda}(\widehat{B})$
and $\ell_{p}^{\lambda}(\widehat{B})$ of non-absolutely  type are
norm isomorphic to the spaces $c_0$, $c$, $\ell_\infty$ and $\ell_p$
, respectively, that is, $c_{0}^{\lambda}(\widehat{B})\cong c_0$,
$c^{\lambda}(\widehat{B})\cong c$,
$\ell_{\infty}^{\lambda}(\widehat{B})\cong\ell_\infty$ and
$\ell_{p}^{\lambda}(\widehat{B})\cong\ell_p$ .
\end{thm}
\begin{pf}
We prove the theorem for the space $c_{0}^{\lambda}(\widehat{B})$.
To prove our assertion we should show the existence of a linear
bijection between the spaces $c_{0}^{\lambda}(\widehat{B})$ and
$c_0$.  Let $T:c_{0}^{\lambda}(\widehat{B})\rightarrow c_0$ be
defined by (\ref{e2.7}) . Then,
$T(x)=y_{k}(\lambda)=\widehat{W}(x)\in c_0$ for every $x\in
c_{0}^{\lambda}(\widehat{B})$ and the linearity of $T$ is clear.
Further, it is trivial that $x=0$ whenever $Tx=\theta$ and hence $T$
is injective .

Moreover, let $y=(y_k)\in c_0$ and define the sequence
$x={x_k(\lambda)}$ by
\begin{equation}\label{e2.6}
x_k(\lambda)=\sum_{j=0}^k\frac{1}{r}\sum_{v=0}^{k-j}\bigg(\frac{-s+\sqrt{s^{2}-4tr}}{2r}\bigg)^{k-j-v}
\bigg(\frac{-s-\sqrt{s^{2}-4tr}}{2r}\bigg)^{v}\sum_{i=j-1}^j(-1)^{j-i}\frac{\lambda_i}{\lambda_j-\lambda_{j-1}}y_i.
\end{equation}
Then we obtain
$$
rx_{k}(\lambda)+sx_{k-1}(\lambda)+tx_{k-2}(\lambda)=\sum_{j=k-1}^k(-1)^{k-j}\frac{\lambda_j}{\lambda_k-\lambda_{k-1}}y_j\qquad
\textrm{for all} \ k\in \mathbb{N}.
$$
Hence, for every $n\in \mathbb{N}$ we get by (\ref{e2.7})
\begin{eqnarray*}
\widehat{W}_n(x)&=&\frac{1}{\lambda_{n}}\sum_{k=0}^{n}
(\lambda_{k}-\lambda_{k-1})(rx_{k}+sx_{k-1}+tx_{k-2})
=\frac{1}{\lambda_n}\sum_{k=0}^n\sum_{j=k-1}^k(-1)^{k-j}\lambda_jy_j=y_n.
\end{eqnarray*}
This show that $\widehat{W}(x)=y$ and since $y\in c_0$, we conclude
that  $\widehat{W}(x)\in c_0$. Thus, we deduce that $x\in
c_0^\lambda(\widehat{B})$ and $Tx=y$. Hence $T$ is surjective.

Moreover one can easily see for every $x\in
c_0^\lambda(\widehat{B})$ that
$$
\|Tx\|_\infty=\|y_{k}(\lambda)\|_\infty=\|\widehat{W}x\|_\infty=\|x\|_{c_{0}^{\lambda}(\widehat{B})}
$$
which means that $T$ is norm preserving. Consequently $T$ is a
linear bijection which show that the spaces
$c_{0}^{\lambda}(\widehat{B})$ and $c_0$ are linearly  isomorphic,
as desired.
\end{pf}

Let $(X,\|.\|)$ be a normed space. A sequence $(b_k)$ of elements of
$X$ is called a Schauder basis for $X$ if and only if, for each
$x\in X$ there exists a unique sequence $(\alpha_k)$ of scalars such
that $x=\sum_k\alpha_kb_k$, i.e. such that
$\lim_{n\rightarrow\infty}\|x-\sum_{k=0}^n\alpha_kb_k\|=0$.

Because of the isomorphism $T$, defined in the proof of Theorem
\ref{t2}, is onto the inverse image of the  basis of those space
$c_{0},c$ and $\ell_{p}$ are the basis of new spaces
$c_{0}^{\lambda}(\widehat{B}),c^{\lambda}(\widehat{B})$ and
$\ell_{p}^{\lambda}(\widehat{B})$ respectively. Therefore, we have
the following:

\begin{thm}
Let $\alpha_k(\lambda)=\widehat{W}_k(x)$ for all $k\in  \mathbb{N}$
and $ \lim_{k\rightarrow\infty}\widehat{W}_k(x)=l$. Define the
sequence $b^{(k)}(\lambda)=\{b_n^{(k)}(\lambda)\}_{n\in\mathbb{N}}$
for every fixed $k\in \mathbb{N}$ by
\begin{equation}\label{e4.3}
b^{(k)}_{n}(\lambda)=\left\{\begin{array}{ll}
  \displaystyle
  \frac{d_{nk}(r,s,t)\lambda_{k}}{\lambda_{k}-\lambda_{k-1}}-\frac{d_{n,k+1}(r,s,t)\lambda_{k}}{\lambda_{k+1}-\lambda_{k}},
   & (n>k)\\\displaystyle\frac{1}{r}\frac{\lambda_{k}}{(\lambda_{k}-\lambda_{k-1})}, & (n=k) \\
   \displaystyle\, 0, & (n<k).
\end{array}\right.
\end{equation}

Then, the following statements hold:\\

(i) The sequence $\{b_n^{(k)}(\lambda)\}_{n\in\mathbb{N}}$ is a
basis for the spaces $c_{0}^{\lambda}(\widehat{B})$ and any $x\in
c_{0}^{\lambda}(\widehat{B})$ has a unique representation
of the norm $x=\sum_k\alpha_k(\lambda)b^{(k)}(\lambda)$.\\

(ii) The sequence $\{b_n^{(k)}(\lambda)\}_{n\in\mathbb{N}}$ is a
basis for the spaces $\ell_{p}^{\lambda}(\widehat{B})$ and any $x\in
\ell_{p}^{\lambda}(\widehat{B})$ has a unique representation
of the norm $x=\sum_k\alpha_k(\lambda)b^{(k)}(\lambda)$.\\

(iii) The sequence $\{b, b^{(0)}(\lambda), b^{(0)}(\lambda),
b^{(1)}(\lambda),... \}$ is a basis for the space
$c^{\lambda}(\widehat{B})$ and any $x\in c^{\lambda}(\widehat{B})$
has a unique representation of the norm
$x=lb+\sum_k[\alpha_k(\lambda)-l]b^{(k)}(\lambda)$, where
$b=(b_k)=\bigg(\sum_{j=0}^kd_{kj}\bigg)$.
\end{thm}

\section{The Inclusion Relations}
In the present section, we prove some inclusion relations concerning
with the spaces
$c_{0}^{\lambda}(\widehat{B}),c^{\lambda}(\widehat{B}),\ell_{\infty}^{\lambda}(\widehat{B})$
and $\ell_{p}^{\lambda}(\widehat{B})$.
\begin{thm}
The inclusion $c_{0}^{\lambda}(\widehat{B})\subset
c^{\lambda}(\widehat{B})$ strictly holds.
\end{thm}
\begin{pf}
It is obvious that the inclusion
$c_{0}^{\lambda}(\widehat{B})\subset c^{\lambda}(\widehat{B})$
holds. Further to show that this inclusion is strict, consider the
sequence $x=(x_k)$ defined by
$$
x_k=\sum_{j=0}^k d_{kj}
$$
for all $k\in \mathbb{N}$. Then
$$
rx_{k}+sx_{k-1}+tx_{k-2}=\sum_{j=k-1}^k(-1)^{k-j}\frac{\lambda_j}{\lambda_k-\lambda_{k-1}}
$$
Then we obtain by (\ref{e2.7}) that
$$
\widehat{W}(x)=\frac{1}{\lambda_{n}}\sum_{k=0}^{n}
(\lambda_{k}-\lambda_{k-1})(rx_{k}+sx_{k-1}+tx_{k-2})
=\frac{1}{\lambda_{n}}\sum_{k=0}^{n} (\lambda_{k}-\lambda_{k-1})=1
$$
for all $n\in \mathbb{N}$, which shows that $\widehat{W}(x)=e$ and
hence $\widehat{W}(x)\in c\backslash c_0$, where $e=(1,1,1,...)$.
Thus ,this sequence $x$ is in $c^{\lambda}(\widehat{B})$ but not in
$c_{0}^{\lambda}(\widehat{B})$. Hence, the inclusion
$c_{0}^{\lambda}(\widehat{B})\subset c^{\lambda}(\widehat{B})$ is
strict and this completes the proof of the theorem.
\end{pf}
\begin{thm}\label{thm5}
If $r+s+t=0$ then the inclusion $c\subset
c_0^{\lambda}(\widehat{B})$ strictly hold.
\end{thm}
\begin{pf}
Suppose that $r+s+t=0$ and $x\in c$. Then
$\widehat{B}x=B(r,s,t)(x)=(rx_{k}+sx_{k-1}+tx_{k-2})\in c_0$ and
hence $\widehat{B}(x)\in c_0^\lambda$ since the inclusion
$c_0\subset c_0^\lambda$ \cite{mmakn1}. This show that $x\in
c_0^\lambda(\widehat{B})$ ,i.e, $c\subset
c_0^{\lambda}(\widehat{B})$ holds. Further consider the sequence
$y=(y_k)$ defined by $y_k=\ln(k+3)$ for all $k\in \mathbb{N}$. Then
it is trivial that $y\notin c$ . On the other hand, it can easily
seen that $\widehat{B}y\in c_0$. Hence $\widehat{B}y\in c_0^\lambda$
which means that $y\in c_0^{\lambda}(\widehat{B})$. Thus the
sequence $y\in c_0^{\lambda}(\widehat{B})\backslash c$. Hence, the
inclusion $c\subset c_0^{\lambda}(\widehat{B})$  is strict. This
completes the proof.
\end{pf}

\begin{lm}\label{lm1}
$A\in (\ell_\infty:c_0)$ if and only if
$\lim_{n\rightarrow\infty}\sum_k|a_{nk}|=0$.
\end{lm}
\begin{thm}\label{thm6}
Let $z=(z_k)$ is defined by
$$
z_k=\bigg|\frac{r(\lambda_{k}-\lambda_{k-1})+s(\lambda_{k+1}-\lambda_{k})+t(\lambda_{k+2}-\lambda_{k+1})}{\lambda_{k}-\lambda_{k-1}}\bigg|
$$
for all $k\in \mathbb{N}$. Then the inclusion $\ell_\infty\subset
c_0^{\lambda}(\widehat{B})$ strictly holds if and only if $z\in
c_0^{\lambda}$.
\end{thm}
\begin{pf}
Let  $\ell_\infty$ be a subset of $c_0^{\lambda}(\widehat{B})$ .
Then we obtain that $\widehat{W}(x)\in c_0$ for every $x\in
\ell_\infty$ and the matrix $\widehat{W}=\{w_{nk}^{\lambda}\}$ is in
the class $(\ell_\infty: c_0)$. By using Lemma \ref{lm1} it follows
that
\begin{equation}\label{e4.0}
\lim_{n\rightarrow \infty}\sum_k|w_{nk}^{\lambda}|=0.
\end{equation}
Now, by taking into account the definition of matrix
$\widehat{W}=\{w_{nk}^{\lambda}\}$ given by (\ref{e2.3}), we have
for every $n\in\mathbb{N}$ that
\begin{eqnarray}\label{e4.2}
\sum_k|w_{nk}^{\lambda}|&=&\frac{1}{\lambda_n}\sum_{k=0}^{n-2}\bigg|r(\lambda_{k}-\lambda_{k-1}) + s(\lambda_{k+1}-\lambda_{k})+t(\lambda_{k+2}-\lambda_{k+1}\bigg|\nonumber\\
&+&\bigg|\frac{r(\lambda_{n-1}-\lambda_{n-2})+s(\lambda_{n}-\lambda_{n-1})}{\lambda_n}\bigg|+|r|\frac{(\lambda_{n}-\lambda_{n-1})}{\lambda_n}
\end{eqnarray}
Thus, the condition (\ref{e4.2}) implies that
\begin{eqnarray}
&&\lim_{n\rightarrow\infty}|r|\frac{(\lambda_{n}-\lambda_{n-1})}{\lambda_n}=0 \label{e4.4}\\
&&\lim_{n\rightarrow\infty}\bigg|\frac{r(\lambda_{n-1}-\lambda_{n-2})+s(\lambda_{n}-\lambda_{n-1})}{\lambda_n}\bigg|=0 \label{e4.5}\\
&&\lim_{n\rightarrow\infty}\frac{1}{\lambda_n}\sum_{k=0}^{n-2}|r(\lambda_{k}-\lambda_{k-1})+s(\lambda_{k+1}-\lambda_{k})+t(\lambda_{k+2}-\lambda_{k+1}|=0.
\label{e4.6}
\end{eqnarray}
Now, we have for every $n\geq 1$ that

$$
\frac{1}{\lambda_n}\sum_{k=0}^{n-2}|r(\lambda_{k}-\lambda_{k-1})+s(\lambda_{k+1}-\lambda_{k})+t(\lambda_{k+2}-\lambda_{k+1}|=\frac{\lambda_{n-2}}{\lambda_n}
\bigg[\frac{1}{\lambda_{n-2}}\sum_{k=0}^{n-2}(\lambda_{k}-\lambda_{k-1})z_k
\bigg]
$$
and since
$\lim_{n\rightarrow\infty}\frac{\lambda_{n-2}}{\lambda_n}=1$ by
(\ref{e4.4}) and (\ref{e4.5}) ; we obtain that by (\ref{e4.6})
\begin{equation}\label{e4.7}
\lim_{n\rightarrow\infty}\frac{1}{\lambda_{n-2}}\sum_{k=0}^{n-2}(\lambda_{k}-\lambda_{k-1})z_k=0
\end{equation}
which shows that $z=(z_k)\in c_0^\lambda$ where the sequence
$z=(z_k)$ is defined by
$$
z_k=\bigg|\frac{r(\lambda_{k}-\lambda_{k-1})+s(\lambda_{k+1}-\lambda_{k})+t(\lambda_{k+2}-\lambda_{k+1})}{\lambda_{k}-\lambda_{k-1}}\bigg|
$$
for all $k\in \mathbb{N}$.

Conversely, we suppose that $z=(z_k)\in c_0^\lambda$. Then we have
(\ref{e4.7}). Further, for every $n\geq 2$, we derive that,
\begin{eqnarray}\label{e4.8}
&&\frac{1}{\lambda_n}\sum_{k=0}^{n-2}|r(\lambda_{k}-\lambda_{k-1})+s(\lambda_{k+1}-\lambda_{k})+t(\lambda_{k+2}-\lambda_{k+1})| \nonumber\\
&=&\frac{1}{\lambda_{n}}\sum_{k=0}^{n-2}(\lambda_{k}-\lambda_{k-1})z_k \nonumber \\
&\leq&\frac{1}{\lambda_{n-2}}\sum_{k=0}^{n-2}(\lambda_{k}-\lambda_{k-1})z_k.
\end{eqnarray}
Then, (\ref{e4.7}) and (\ref{e4.8}) together imply that (\ref{e4.6})
holds. On the other hand, we have for every $n\geq 2$ that

\begin{eqnarray*}
&&\bigg|\frac{r\lambda_{n-2}+s(\lambda_{n-1}-\lambda_{0})+t(\lambda_{n}-\lambda_{1})}{\lambda_n}\bigg|\\
&=&\bigg|\frac{1}{\lambda_n}\sum_{k=0}^{n-2}r(\lambda_{k}-\lambda_{k-1})+s(\lambda_{k+1}-\lambda_{k})
+t(\lambda_{k+2}-\lambda_{k+1})\bigg| \\
&\leq&\frac{1}{\lambda_n}\sum_{k=0}^{n-2}|r(\lambda_{k}-\lambda_{k-1})+s(\lambda_{k+1}-\lambda_{k})+t(\lambda_{k+2}-\lambda_{k+1})|
\end{eqnarray*}
Therefore, it follows by (\ref{e4.6})
$$\lim_{n\rightarrow\infty}\frac{r\lambda_{n-2}+s(\lambda_{n-1}-\lambda_{0})+t(\lambda_{n}-\lambda_{1})}{\lambda_n}=0.$$
Particularly if (i) $r=0$, $s=-1$, $t=1$, then we obtain
$\lim_{n\rightarrow\infty}[(\lambda_{n}-\lambda_{n-1})/\lambda_n]=0$,
which shows that (\ref{e4.4}) holds. (ii) $r=-1$, $s=1$, $t=0$, then
we obtain
$\lim_{n\rightarrow\infty}[(\lambda_{n-1}-\lambda_{n-2})/\lambda_n]=0$,
which shows that (\ref{e4.5}) holds. Thus, we deduce by the relation
(\ref{e4.2}) that (\ref{e4.0}) holds. This leads us with
Lemma\ref{lm1} to the consequence that,
$\widehat{W}(x)\in(\ell_\infty: c_0)$. Hence, the inclusion
$\ell_\infty\subset c_0^{\lambda}(\widehat{B})$ holds. Lastly, it is
obvious that the sequence $y$, defined in the proof of
Theorem\ref{thm5}, is in $c_0^{\lambda}(\widehat{B})$ but not in
$\ell_\infty$, so the inclusion $\ell_\infty\subset
c_0^{\lambda}(\widehat{B})$ is strict.

\end{pf}

\begin{thm}
Let $1\leq p<\infty$. Then the inclusions
$\ell_\infty\subset\ell_{\infty}^{\lambda}(\widehat{B})$ and
$\ell_{p}^{\lambda}(\widehat{B})\subset\ell_{\infty}^{\lambda}(\widehat{B})$
are strictly holds.
\end{thm}
\begin{pf}
To prove the validity of the inclusion
$\ell_\infty\subset\ell_{\infty}^{\lambda}(\widehat{B})$, it
suffices to show that, for every $x\in
\ell_{\infty}^{\lambda}(\widehat{B})$, there exists a positive real
number $K$ such that
$\|x\|_{\ell_{\infty}^{\lambda}(\widehat{B})}\leq
K\|x\|_{\ell_\infty}$. Let $x\in \ell_{\infty}$. Then, we have

\begin{eqnarray*}
\|x\|_{\ell_{\infty}^{\lambda}(\widehat{B})}&=&\sup_{n\in\mathbb{N}}\bigg|\frac{1}{\lambda_{n}}\sum_{k=0}^{n}
(\lambda_{k}-\lambda_{k-1})(rx_{k}+sx_{k-1}+tx_{k-2})\bigg| \\
&\leq&\sup_{0\leq k\leq
n}|rx_{k}+sx_{k-1}+tx_{k-2}|\sup_{n\in\mathbb{N}}\bigg|\frac{1}{\lambda_{n}}\sum_{k=0}^{n}
(\lambda_{k}-\lambda_{k-1})\bigg|\\
&\leq& \max \{|r|+|s|+|t|\}\|x\|_{\ell_\infty}
\end{eqnarray*}
so that $x\in \ell_{\infty}^{\lambda}(\widehat{B})$ and hence,
$\ell_\infty\subset\ell_{\infty}^{\lambda}(\widehat{B})$.
Furthermore, we consider $x=(x_k)$ defined by
$x_k=\sum_{j=0}^kd_{kj}$ for all $k\in \mathbb{N}$. Then we have
$\widehat{W}(x)=1$. Thus, we deduce that
$\|x\|_{\ell_{\infty}^{\lambda}(\widehat{B})}=
\|\widehat{W}x\|_{\infty}=1$. Hence $x=(x_k)\in
\ell_{\infty}^{\lambda}(\widehat{B})$.

On the other hand, we know that the inclusion $\ell_\infty\subset
c_0^{\lambda}(\widehat{B})$ is strict from Theorem \ref{thm6}. Since
$c_0^{\lambda}(\widehat{B})\subset\ell_\infty(\widehat{B})$, the
inclusion $\ell_\infty\subset\ell_\infty(\widehat{B})$ strictly
holds.

 Now, we suppose that $x\in \ell_p^\lambda(\widehat{B})$. Since $\widehat{W}(x)\in
 \ell_p\subset\ell_\infty$, we have $x\in \ell_\infty^\lambda(\widehat{B})$.
 Hence, $\ell_{p}^{\lambda}(\widehat{B})\subset\ell_{\infty}^{\lambda}(\widehat{B})$.
 Further, to show that this inclusion is strict, we consider the
 sequence $x=(x_k)$ defined by $x=(1,1,1,...)$ and assume that
 $r+s+t=1$. Then, we have

 \begin{eqnarray*}
\|\widehat{W}x\|_{\infty}&=
\sup_{n\in\mathbb{N}}\bigg|\frac{1}{\lambda_{n}}\sum_{k=0}^{n}
(\lambda_{k}-\lambda_{k-1})(rx_{k}+sx_{k-1}+tx_{k-2})\bigg| \\&=1.
 \end{eqnarray*}
But, for $1\leq p<\infty$
$$
\sum_{n=0}^\infty\bigg|\frac{1}{\lambda_{n}}\sum_{k=0}^{n}
(\lambda_{k}-\lambda_{k-1})(rx_{k}+sx_{k-1}+tx_{k-2})\bigg|^p=\sum_{n=0}^\infty1=\infty.
$$
Hence, we have $x\notin \ell_p^\lambda(\widehat{B})$, which we want
to show.

\end{pf}

\section{The $\alpha$-, $\beta$- and $\gamma$-Duals of The Spaces
$c_{0}^{\lambda}(\widehat{B}), c^{\lambda}(\widehat{B}),
\ell_{\infty}^{\lambda}(\widehat{B})$ and
$\ell_{p}^{\lambda}(\widehat{B})$}
In this section, we determine the $\alpha$-, $\beta$- and
$\gamma$-duals of the generalized difference sequence spaces
$c_{0}^{\lambda}(\widehat{B}), c^{\lambda}(\widehat{B}),
\ell_{\infty}^{\lambda}(\widehat{B})$ and
$\ell_{p}^{\lambda}(\widehat{B})$ of non-absolute type.

We shall firstly give the definition of $\alpha$-, $\beta$- and
$\gamma$-duals of a sequences space and after quote the lemmas due
to Stieglitz and Tietz \cite{msht} which are needed in proving the
theorems given in Section 4 and 5.

For the sequence spaces $\lambda$ and $\mu$, define the set
$S(\lambda,\mu)$ by
\begin{equation}\label{e5.1}
S(\lambda,\mu)=\{z=(z_k)\in w: xz=(x_kz_k)\in \mu \ \textrm{for
all}\ x\in \lambda\}.
\end{equation}

With the notation of (\ref{e5.1}), the $\alpha$-, $\beta$- and
$\gamma$-duals of a sequences space $\lambda$, which are
respectively denoted by
$$
\lambda^\alpha=S(\lambda,\ell_1),\quad
\lambda^\beta=S(\lambda,cs)\quad \textrm{and}\ \
\lambda^\gamma=S(\lambda,bs).
$$

\begin{lm}\label{e4.001}
(i) $A\in (c_0:\ell_1)=(c:\ell_1)=(\ell_\infty:\ell_1)$ if and only
if
$$
 \sup_{K\in \mathcal{F}}\sum_n\bigg|\sum_{k\in
K}a_{nk}\bigg|<\infty.
$$

(ii) Let $1<p<\infty$ and $\frac{1}{p}+\frac{1}{q}=1$. Then  $A\in
(\ell_p:\ell_1)$ if and only if
$$
 \sup_{N\in \mathcal{F}}\sum_k\bigg|\sum_{n\in
N}a_{nk}\bigg|^{q}<\infty;
$$
\end{lm}

\begin{lm}\label{e4.10}
$A\in (c_0:c)$ if and only if
\begin{equation}\label{e5.2}
\lim_{n\rightarrow\infty} a_{nk}=\alpha_k \ \textrm{for each fixed}\
k\in \mathbb{N},
\end{equation}

\begin{equation}\label{e5.3}
\sup_{n\in \mathbb{N}}\sum_k|a_{nk}|<\infty.
\end{equation}
\end{lm}

\begin{lm}\label{e4..2}
$A\in (c:c)$ if and only if (\ref{e5.2}) and (\ref{e5.3}) hold, and
\begin{equation}
\lim_{n\rightarrow\infty}\sum_ka_{nk}\quad\ \textrm{exists}.
\end{equation}
\end{lm}

\begin{lm}
$A\in (\ell_\infty:c)$ if an only if  (\ref{e5.2}) holds and
$$
\lim_{n\rightarrow\infty}\sum_k|a_{nk}|=\sum_k|\alpha_k|.
$$
\end{lm}

\begin{lm}\label{e4.003}
Let  $1<p<\infty$. Then, $A\in (\ell_p:c)$ if and only if
(\ref{e5.2}) hold and
\begin{equation}\label{e4..4}
\sup_{n\in \mathbb{N}}\sum_k|a_{nk}|^q<\infty\qquad
\bigg(\frac{1}{p}+\frac{1}{q}=1\bigg).
\end{equation}
\end{lm}

\begin{lm}\label{e4..1}
$A\in (c:\ell_\infty)=(c_0:\ell_\infty)=(\ell_\infty:\ell_\infty)$
if and only if (\ref{e5.3}) holds.
\end{lm}

\begin{lm}\label{e4..3}
Let  $1<p<\infty$. Then, $A\in (\ell_p:\ell_\infty)$ if and only if
(\ref{e4..4}) holds.

\end{lm}

Now we consider  the following sets:

\begin{eqnarray*}
f_1^\lambda&=&\bigg\{ a=(a_n)\in w:\sup_{K\in
\mathcal{F}}\sum_n\bigg|\sum_{k\in
K}f_{nk}^\lambda\bigg|<\infty\bigg\},\\
f_2^\lambda&=&\bigg\{ a=(a_n)\in w:\sum_{j=k}^\infty d_{jk}a_j \
\textrm{exists for
each}\ k\in \mathbb{N}.\bigg\},\\
 f_3^\lambda&=&\bigg\{
a=(a_n)\in w:\sup_{n\in \mathbb{N}}\sum_{k=0}^{n-1}|\widehat{g}_k(n)|^{q}<\infty\bigg\},\\
f_4^\lambda&=&\bigg\{ a=(a_n)\in
w:\sup_{n\in\mathbb{N}}\bigg|\frac{1}{r}\frac{\lambda_n}{(\lambda_n-\lambda_{n-1})}a_n\bigg|<\infty\bigg\},\\
f_5^\lambda&=&\bigg\{ a=(a_n)\in
w:\lim_{n\rightarrow\infty}\sum_{k=0}^n\bigg[\sum_{j=0}^k
d_{kj}\bigg]a_k \quad \textrm{exists}\},\\
f_6^\lambda&=&\bigg\{ a=(a_n)\in w:\sup_{N\in
\mathcal{F}}\sum_{k=0}^\infty\bigg|\sum_{n\in
N}f_{nk}^\lambda\bigg|^q<\infty\bigg\},
\\
f_7^\lambda&=&\bigg\{ a=(a_n)\in w:\lim_{n\rightarrow \infty}
\sum_{k} |v_{nk}^{\lambda}|=\sum_{k} |\lim_{n\rightarrow \infty}
v_{nk}^{\lambda}|\bigg\},
\end{eqnarray*}
where the matrices $F^\lambda=(f^{\lambda}_{nk})$ and
$V^\lambda=(v^\lambda_{nk})$ are defined as follows,
$$
f^{\lambda}_{nk}=\left\{\begin{array}{ll}
  \displaystyle\bigg(\frac{d_{nk}\lambda_{k}}{\lambda_{k}-\lambda_{k-1}}-\frac{d_{n,k+1}\lambda_{k}}{\lambda_{k+1}-\lambda_{k}}\bigg)a_{n},
   & (k<n)\\  \displaystyle\frac{1}{r}\frac{\lambda_{n}}{(\lambda_{n}-\lambda_{n-1})}a_{n}, & (k=n) \\\displaystyle\ 0, & (k>n)
\end{array}\right.
$$
 and

$$
v^{\lambda}_{nk}=\left\{\begin{array}{ll}
  \displaystyle\widehat{g}_k(n),
   & (k<n)\\  \displaystyle\frac{1}{r}\frac{\lambda_{n}}{(\lambda_{n}-\lambda_{n-1})}a_n, & (k=n) \\\displaystyle\ 0, & (k>n)
\end{array}\right.
$$

for all $k,n\in \mathbb{N}$ and the $\widehat{g}_k(n)$ is defined as
follows
$$
\widehat{g}_k(n)=\lambda_k\bigg(\frac{1}{\lambda_k-\lambda_{k-1}}\sum_{j=k}^n
d_{jk}a_j- \frac{1}{\lambda_{k+1}-\lambda_{k}}\sum_{j=k+1}^n
d_{j,k+1}a_j\bigg)\qquad \textrm{for}\ k<n.
$$

\begin{thm}\label{e4.002}
(i)
$\{c_{0}^{\lambda}(\widehat{B})\}^\alpha=\{c^{\lambda}(\widehat{B})\}^\alpha
=\{\ell_\infty^\lambda(\widehat{B})\}^\alpha=f_1^\lambda$.

(ii) Let $1<p<\infty$ and $\frac{1}{p}+\frac{1}{q}=1$. Then,
$\{\ell_p^\lambda(\widehat{B})\}^\alpha=f_6^\lambda$.
\end{thm}
\begin{pf}
We prove the theorem for the space $c_{0}^{\lambda}(\widehat{B})$.
Let $a=(a_n)\in w$. Then, we obtain the equality
\begin{equation}\label{e5.4}
a_nx_n=\sum_{k=0}^nd_{nk}\sum_{j=k-1}^k(-1)^{k-j}\frac{\lambda_j}{\lambda_k-\lambda_{k-1}}a_ny_i=F_n^\lambda(y)\qquad
\textrm{for all}\ n\in \mathbb{N},
\end{equation}
 by relation (\ref{e2.6}). Thus we observe by
 (\ref{e5.4}) that $ax=(a_nx_n)\in\ell_1$ whenever $x=(x_k)\in
 c_{0}^{\lambda}(\widehat{B})$ if and only if $F^\lambda y\in \ell_1$
 whenever $y=(y_k)\in c_0$ . This means that the sequence $a=(a_n)\in \{c_{0}^{\lambda}(\widehat{B})\}^\alpha$
 if and only if
 $F^\lambda\in(c_0:\ell_1)$. Therefore we  obtain by
 Lemma \ref{e4.001} with $F^\lambda$ instead of $A$ that $a=(a_n)\in \{c_{0}^{\lambda}(\widehat{B})\}^\alpha$ if and only if
 $$
 \sup_{K\in
\mathcal{F}}\sum_n\bigg|\sum_{k\in K}f_{nk}^\lambda\bigg|<\infty
 $$
 which leads us to the consequence that $ \{c_{0}^{\lambda}(\widehat{B})\}^\alpha=
f_1^\lambda$. This completes
 the proof.
\end{pf}

\begin{thm}\label{e4.004}
(i) $\{c_{0}^{\lambda}(\widehat{B})\}^\beta= f_2^\lambda\cap
f_3^\lambda\cap f_4^\lambda \quad (\textrm{with} \quad q=1)$.

(ii) $\{c^{\lambda}(\widehat{B})\}^\beta=f_2^\lambda\cap
f_3^\lambda\cap f_4^\lambda\cap f_5^\lambda \quad (\textrm{with}
\quad q=1)$.

(iii)
$\{\ell_{\infty}^{\lambda}(\widehat{B})\}^{\beta}=f_2^\lambda\cap
f_4^\lambda\cap f_7^\lambda$

(iv) Let $1<p<\infty$ and $\frac{1}{p}+\frac{1}{q}=1$. Then,
$\{\ell_p^\lambda(\widehat{B})\}^\beta=f_2^\lambda\cap f_3^\lambda
\cap f_4^\lambda$.
\end{thm}

\begin{pf}
Consider the equality
\begin{eqnarray}
\sum_{k=0}^{n}a_kx_k&=&\sum_{k=0}^{n}\bigg[\sum_{j=0}^{k}d_{kj}
\sum_{i=j-1}^{i}(-1)^{j-i}\frac{\lambda_{i}}{\lambda_{j}-\lambda_{j-1}}y_i\bigg]a_k \nonumber \\
\quad\quad\quad\quad\quad\quad\quad\quad\quad\quad\quad &=&
\sum_{k=0}^{n}(\lambda_{k}y_{k}-\lambda_{k-1}y_{k-1})\bigg[\frac{1}{\lambda_{k}-\lambda_{k-1}}\sum_{j=k}^{n}d_{jk}a_j\bigg] \nonumber\\
&=&\sum_{k=0}^{n-1}\lambda_k\bigg[\frac{\sum_{j=k}^{n}d_{jk}a_j}{\lambda_{k}-\lambda_{k-1}}-\frac{\sum_{j=k+1}^{n}d_{j,k+1}a_j}{\lambda_{k+1}-\lambda_{k}}\bigg]
y_k+\frac{1}{r}\frac{a_n\lambda_n}{(\lambda_{n}-\lambda_{n-1})}y_n \nonumber\\
&=&\sum_{k=0}^{n-1}\widehat{g}_k(n)y_k+\frac{1}{r}\frac{a_n\lambda_n}{(\lambda_{n}-\lambda_{n-1})}y_n \nonumber\\
&=&V^\lambda_n(y) ; \quad (n\in\mathbb{N})\label{e4.9}.
\end{eqnarray}
 Then we deduce by (\ref{e4.9}) that
$ax=(a_kx_k)\in cs$ whenever $x=(x_k)\in c_0^\lambda(\widehat{B})$
if and only if $V^\lambda y\in c$ whenever $y=(y_k)\in c_0$. This
means that $a=(a_k)\in \{c_0^\lambda(\widehat{B})\}^\beta$ if and
only if $V^\lambda\in (c_0:c)$. Therefore, by using Lemma
\ref{e4.10}, we obtain :
\begin{eqnarray}
&&\sum_{j=k}^{\infty}d_{jk}a_j\qquad \textrm{exists for each}\
k\in\mathbb{N},\\
&&\sup_{n\in\mathbb{N}}\sum_{k=0}^{n-1}|\widehat{g}_k(n)|<\infty,\\
&&\sup_{k\in\mathbb{N}}\bigg|\frac{1}{r}\frac{\lambda_{k}}{(\lambda_{k}-\lambda_{k-1})}a_k\bigg|<\infty.
\end{eqnarray}
Hence, we conclude that
 $\{c_0^\lambda(\widehat{B})\}^\beta=f_2^\lambda\cap
f_3^\lambda\cap f_4^\lambda$.
\end{pf}

Finally, we ended up this section with the following theorem which
determines the $\gamma$-duals of sequence spaces
$c_{0}^{\lambda}(\widehat{B}), c^{\lambda}(\widehat{B}),
\ell_{\infty}^{\lambda}(\widehat{B})$ and
$\ell_{p}^{\lambda}(\widehat{B})$ .

\begin{thm}\label{e4.005}
(i) $ \{c_{0}^{\lambda}(\widehat{B})\}^\gamma=
 \{c^{\lambda}(\widehat{B})\}^\gamma=\{\ell_\infty^{\lambda}(\widehat{B})\}^\gamma=f_3^\lambda\cap f_4^\lambda
 \quad (\textrm{with} \quad q=1)$.

 (ii)  Let $1<p<\infty$ and $\frac{1}{p}+\frac{1}{q}=1$. Then,
 $\{\ell_p^\lambda(\widehat{B})\}^\gamma=f_3^\lambda \cap
 f_4^\lambda$.
\end{thm}

\section{Some Matrix Transformations Related To The Spaces
$c_{0}^{\lambda}(\widehat{B}), c^{\lambda}(\widehat{B})$,
$\ell_{\infty}^{\lambda}(\widehat{B})$ and
$\ell_{p}^{\lambda}(\widehat{B})$}

In this final section, we state some results which characterize
various matrix mappings on the spaces $c_{0}^{\lambda}(\widehat{B}),
c^{\lambda}(\widehat{B}), \ell_{\infty}^{\lambda}(\widehat{B})$ and
$\ell_{p}^{\lambda}(\widehat{B})$. We shall write throughout for
brevity that

$$
\widehat{g}_{nk}(m)=\lambda_k\bigg(\frac{1}{\lambda_k-\lambda_{k-1}}\sum_{j=k}^m
d_{jk}a_{nj}- \frac{1}{\lambda_{k+1}-\lambda_{k}}\sum_{j=k+1}^m
d_{j,k+1}a_{nj}\bigg)\qquad \textrm{for}\ k<m.
$$
and

$$
\widehat{g}_{nk}=\lambda_k\bigg(\frac{1}{\lambda_k-\lambda_{k-1}}\sum_{j=k}^\infty
d_{jk}a_{nj}- \frac{1}{\lambda_{k+1}-\lambda_{k}}\sum_{j=k+1}^\infty
d_{j,k+1}a_{nj}\bigg)
$$
for all $k,m,n\in \mathbb{N}$ provided the convergence of the
series.

We shall begin with lemmas which are needed in the proof of our
theorems.

\begin{lm}\label{e5..1}
The matrix mappings between the $BK-$spaces are continuous.
\end{lm}
\begin{lm}\label{e5.5}
$A=(a_{nk})\in (c:\ell_p)$ if and only if
\begin{equation}
\sup_{F\in \mathcal{F}} \sum_{n}\bigg|\sum_{k\in F}
a_{nk}\bigg|^p<\infty; \quad (1<p<\infty)
\end{equation}

\end{lm}
\begin{lm}\label{lm11}
$A=(a_{nk})\in (c:c_0)$ if and only if (\ref{e5.3}) holds and
\begin{eqnarray}\label{e5-2}
\quad\quad\quad\quad\quad\quad\quad\quad\quad\quad\lim_{n\rightarrow\infty}
a_{nk}=0\quad\quad\ \textrm{for each fixed}\ k\in \mathbb{N}
\end{eqnarray}
\begin{eqnarray*}
\lim_{n\rightarrow\infty}\sum_k a_{nk}=0
\end{eqnarray*}
\end{lm}
\begin{lm}\label{lm12}
$A\in (c_0:c_0)$ if and only if (\ref{e5.3}) and (\ref{e5-2}) hold
\end{lm}

Now, we may give our matrix transformations.
\begin{thm}\label{e5..0}
Let $A=(a_{nk})$ be an infinite matrix and $1<p<\infty$ . Then,
$A\in (c^\lambda(\widehat{B}):\ell_p)$ if and only if
\begin{eqnarray}
\label{e5.7}&&\sup_{F\in \mathcal{F}} \sum_{n}\bigg|\sum_{k\in F}
\widehat{g}_{nk}\bigg|^p<\infty; \quad (1<p<\infty)\\
\label{e5.8}&&\sum_{j=k}^\infty d_{jk}a_{nj}\quad\quad
\textrm{exists for all
fixed}\ k,n\in \mathbb{N}\\
\label{e5.9}&& \sup_{m\in \mathbb{N}} \sum_{k=0}^{m-1}
|\widehat{g}_{nk}(m)|<\infty\\
\label{e5.10}&&\lim_{k\rightarrow\infty}\frac{1}{r}\frac{\lambda_k}{(\lambda_{k}-\lambda_{k-1})}a_{nk}=a_n
\quad\quad \textrm{exists for each fixed}\ n\in \mathbb{N}\\
\label{e5.11}&&\bigg\{\sum_{j=k}^\infty d_{jk}a_{nj}\bigg\}\in cs; \quad (n\in \mathbb{N})\\
\label{e5.12}&& (a_n)\in\ell_p
\end{eqnarray}
\end{thm}
\begin{pf}
If  the conditions (\ref{e5.7})-(\ref{e5.12}) hold and $x=(x_k)$ be
any sequence in the space $c^\lambda(\widehat{B})$ then by using
Theorem \ref{e4.004} we have  that $\{a_{nk}\}_{k\in\mathbb{N}}\in
\{c^\lambda(\widehat{B})\}^\beta$ for all $n\in \mathbb{N}$. Hence,
$A-$transform of $x$ exists,i.e., $Ax$ exists. Furthermore, since
the associated sequence $y=(y_k)$ is in the space $c$, we may write
$\lim_k y_k=l$ for some suitable $l$. Also, the matrix
$\widehat{A}=(\widehat{g}_{nk})$ is in the class $(c:\ell_p)$ by
Lemma \ref{e5.5} and condition (\ref{e5.7}) where $1<p<\infty$.

Now, we may consider the $m^{th}$ partial sum of the series $\sum_k
 a_{nk}x_k$ which is derived by using the relation (\ref{e2.6}):
\begin{equation}\label{e5.13}
\sum_{k=0}^{m}a_{nk}x_k=\sum_{k=0}^{m-1}\widehat{g}_{nk}(m)y_k
+\frac{1}{r}\frac{\lambda_m}{(\lambda_{m}-\lambda_{m-1})}a_{nm}y_m;\quad
\textrm{for all} \quad n,m\in \mathbb{N}
\end{equation}
Then, since $y\in c$ and $\widehat{A}\in (c:\ell_p)$; $\widehat{A}y$
exists and so the series  $\sum_k \widehat{g}_{nk}y_k$ converges for
every $n\in \mathbb{N}$. Also, it follows by (\ref{e5.8}) we have
that
$\lim_{m\rightarrow\infty}\widehat{g}_{nk}(m)=\widehat{g}_{nk}$.
Therefore, if we pass to limit in (\ref{e5.13}) as
$m\rightarrow\infty$, then we obtain by (\ref{e5.10}) that
\begin{equation}\label{e5.14}
\sum_{k}a_{nk}x_k=\sum_{k}\widehat{g}_{nk}y_k +la_n\quad \textrm{for
all} \quad n\in \mathbb{N}
\end{equation}
which can be written as follows :
\begin{equation}\label{e5.15}
A_n(x)=\widehat{A}_n(y) + la_n \quad \textrm{for all} \quad n\in
\mathbb{N}.
\end{equation}
This yields by taking $\ell_p-$norm that
$$
\|Ax\|_{\ell_p}\leq
\|\widehat{A}y\|_{\ell_{p}}+|l|\|a_n\|_{\ell_{p}}<\infty.
$$
Consequently, we have that $Ax\in \ell_p$, i.e., $A\in
(c^\lambda(\widehat{B}):\ell_p)$.

Conversely, if $A\in (c^\lambda(\widehat{B}):\ell_p)$, we have
$\{a_{nk}\}_{k\in N}\in \{c^\lambda(\widehat{B})\}^\beta$ for all
$n\in\mathbb{N}$ which implies with Theorem\ref{e4.004} that the
conditions (\ref{e5.8}), (\ref{e5.9}) and (\ref{e5.11}) are
necessary.

On the other hand, since $c^\lambda(\widehat{B})$ and $\ell_p$ are
$BK-$spaces, we have by Lemma \ref{e5..1} that there is a constant
$M>0$ such that

\begin{equation}\label{e5..2}
\|Ax\|_p\leq M\|x\|_{c^\lambda(\widehat{B})}
\end{equation}
holds for all $x\in c^\lambda(\widehat{B})$. Now $F\in \mathcal{F}$.
Then, the sequence $z=\sum_{k\in F}b^{(k)}(\lambda)$ is in $
c^\lambda(\widehat{B})$, where the sequences
$b^{(k)}(\lambda)=\{b_n^{(k)}(\lambda)\}_{n\in\mathbb{N}}$ is
defined by (\ref{e4.3}) for every fixed $k\in \mathbb{N}$.

Since $\widehat{W}(b^{(k)}(\lambda))=e^{(k)}$ for each fixed $k\in
\mathbb{N}$, we have
$$
\|z\|_{c^\lambda(\widehat{B})}=\|\widehat{W}(z)\|_{\ell_{\infty}}
=\bigg\|\sum_{k\in
F}\widehat{W}(b^{(k)}(\lambda))\bigg\|_{\ell_{\infty}}=\bigg\|\sum_{k\in
F}e^{(k)}\bigg\|_{\ell_{\infty}}=1
$$

Furthermore, for every $n\in \mathbb{N}$, we obtain by (\ref{e4.3})
that
$$
A_n(z)=\sum_{k\in F}A_n(b^{(k)}(\lambda))=\sum_{k\in
F}\sum_ja_{nj}b_j^{(k)}(\lambda)=\sum_{k\in F}\widehat{g}_{nk}
$$
Hence, since the inequality (\ref{e5..2}) is satisfied for the
sequence $z\in c^\lambda(\widehat{B})$ ; we have for any $F\in
\mathcal{F}$ that
$$
\bigg(\sum_n\bigg|\sum_{k\in
F}\widehat{g}_{nk}\bigg|^p\bigg)^{1/p}\leq M
$$
which shows the necessity of (\ref{e5.7}). Thus, it follows by
Lemma\ref{e5.5} that $\widehat{A}=(\widehat{g}_{nk})\in (c:\ell_p)$.

Moreover, we consider the sequence $x=x_k$ defined by (\ref{e2.6})
for every $k\in \mathbb{N}$ and suppose that $y=(y_k)\in c\backslash
c_0$. Then, since $x\in c^\lambda(\widehat{B})$ such that
$y=\widehat{W}(x)$ by (\ref{e1..1}), the transforms $Ax$
and$\widehat{A}y$ exists. Hence, the series $\sum_k a_{nk}x_k$ and
$\sum\widehat{g}_{nk}y_k$ converges for every $n\in \mathbb{N}$. So
we infer that
$$
\lim_{m\rightarrow\infty}\sum_{k=0}^{m-1}\widehat{g}_{nk}(m)y_k=
\sum\widehat{g}_{nk}y_k;\quad (n\in\mathbb{N}).
$$
Consequently, we obtain from (\ref{e5.13}) that
$$
\lim_{m\rightarrow\infty}\frac{1}{r}\frac{\lambda_m}{(\lambda_{m}-\lambda_{m-1})}a_{nm}y_m\quad
\textrm{exists for each fixed} \quad n\in \mathbb{N}.
$$
Hence, we deduce that
$$
\lim_{m\rightarrow\infty}\frac{1}{r}\frac{\lambda_m}{(\lambda_{m}-\lambda_{m-1})}a_{nm}\quad
\textrm{exists for each fixed} \quad n\in \mathbb{N}.
$$
which leads us to the necessity of (\ref{e5.10}) and so the relation
(\ref{e5.13}) holds, where $l=\lim_k y_k$.

Finally, since   $Ax\in \ell_p$ and $\widehat{A}y\in \ell_p$; the
necessity of (\ref{e5.12}) is immediate by (\ref{e5.13}). This
completes the proof.
\end{pf}

\begin{thm}\label{e5...4}
In order that  $A=(a_{nk})\in (c^\lambda(\widehat{B}):\ell_\infty)$
where $A=(a_{nk})$ be an infinite matrix , it is necessary and
sufficient that  (\ref{e5.10}) and (\ref{e5.11}) hold, and
\begin{eqnarray}\label{e5..5}
\label{e5..3}&&\sup_{n\in\mathbb{N}}\sum_k|\widehat{g}_{nk}|<\infty,\\
\label{e5..4}&&(a_n)\in \ell_\infty.
\end{eqnarray}
\end{thm}

\begin{pf}
This is an immediate consequence of Lemma \ref{e4..1} and Theorem
\ref{e5..0}.
\end{pf}

\begin{thm}\label{th13}
(i)In order that $A=(a_{nk})\in (c_0^\lambda(\widehat{B}):\ell_p)$,
it is necessary and sufficient that (\ref{e5.7}),(\ref{e5.8}) and
(\ref{e5.9}) hold and
\begin{equation}\label{e5..6}
\sup_k\bigg|\frac{1}{r}\frac{\lambda_k}{(\lambda_{k}-\lambda_{k-1})}a_{nk}\bigg|<\infty
 \quad\quad (n\in \mathbb{N}).
\end{equation}
(ii))In order that $A\in (c_0^\lambda(\widehat{B}):\ell_\infty)$, it
is necessary and sufficient that (\ref{e5.9}), (\ref{e5..5}) and
(\ref{e5..6}) hold.
\end{thm}

\begin{pf}
This may be obtained by proceedings as in Theorem \ref{e5..0} ,
above. So, we omit the detail.
\end{pf}
\begin{thm}\label{e5...5}
In order that $A=(a_{nk})\in (c^\lambda(\widehat{B}):c)$, it is
necessary and sufficient that  (\ref{e5.10}),
(\ref{e5.11}), (\ref{e5..5})  hold and
\begin{eqnarray}
\label{e5...0}&&\lim_na_n=a\\
\label{e5...1}&&\lim_n\widehat{g}_{nk}=\alpha_k\quad\quad
 \textrm{ for each fixed}\ k\in \mathbb{N}\\
\label{e5...2}&&\lim_n\sum_k\widehat{g}_{nk}=\alpha
\end{eqnarray}
\end{thm}
\begin{pf}
Suppose that $A$ satisfies conditions (\ref{e5.10}), (\ref{e5.11}),
(\ref{e5..5}),(\ref{e5...0}), (\ref{e5...1}) and  (\ref{e5...2}),
and let $x=(x_k)$ be a sequence in the space
$c^\lambda(\widehat{B}$). Since (\ref{e5..5}) implies (\ref{e5.9});
we have by Theorem\ref{e4.004} that $\{a_{nk}\}_{k\in\mathbb{N}}\in
\{c^\lambda(\widehat{B})\}^\beta$ for all $n\in \mathbb{N}$. Hence,
$Ax$ exists . We also observe from (\ref{e5..5}) and (\ref{e5...1})
that
$$
\sum_{j=0}^k|\alpha_j|\leq\sup_{n\in
\mathbb{N}}\sum_j|\widehat{g}_{nk}|<\infty
$$
holds for every $k\in\mathbb{N}$. So $(\alpha_k)\in\ell_1$ and hence
the series $\sum_k\alpha_k(y_k-l)$ converges, where $y=(y_k)\in c$
is the sequence connected with $x=(x_k)$ by the relation
(\ref{e1..1}) such that $\lim_k y_k=l$. Also, it is obvious by
combining  Lemma \ref{e4..2} with the conditions (\ref{e5..5}),
(\ref{e5...1}) and (\ref{e5...2}) that the matrix
$\widehat{A}=(\widehat{g}_{nk})$ is in the class $(c:c)$.

Now, by following the similar way used in the proof of
Theorem\ref{e5..0}, we obtain that the relation (\ref{e5.14}) holds,
which can be written as follows:
\begin{equation}\label{e5...3}
\sum_{k}a_{nk}x_k=\sum_{k}\widehat{g}_{nk}(y_k-l)
+l\sum_{k}\widehat{g}_{nk}+la_n\quad \textrm{for all} \quad n\in
\mathbb{N}
\end{equation}
If we pass the limit in (\ref{e5...3}) as $n\rightarrow \infty$ we
have that
$$
\lim_nA_n(x)=\sum_{k}\alpha_k(y_k-l)+l(\alpha+a),
$$
which shows that $Ax\in c$, i.e.,  $A\in
(c^\lambda(\widehat{B}):c)$.

Conversely, suppose that  $A\in (c^\lambda(\widehat{B}):c)$. Since
the inclusion $c\subset\ell_\infty$ holds; $A\in
(c^\lambda(\widehat{B}):\ell_\infty)$. Therefore, the necessity of
conditions (\ref{e5.10}), (\ref{e5.11}) and  (\ref{e5..5}) are
obvious from Theorem\ref{e5...4}. Furthermore, consider the sequence
$b^{(k)}(\lambda)=\{b_n^{(k)}(\lambda)\}_{n\in\mathbb{N}}\in
c^\lambda(\widehat{B})$  defined by (\ref{e4.3}) for every fixed
$k\in \mathbb{N}$.  Then, one can see  that
$Ab^{(k)}(\lambda)=\{\widehat{g}_{nk}\}_{n\in\mathbb{N}}$ and hence
$\{\widehat{g}_{nk}\}_{n\in\mathbb{N}}\in c$ for every
$k\in\mathbb{N}$ which shows that the necessity of (\ref{e5...1}).
Let $z=\sum_k b^{(k)}(\lambda)$. Then , since the linear
transformation $T:c^\lambda(\widehat{B})\rightarrow c$, defined as
in the proof of Theorem\ref{t2} by analogy, is continuous and
$\widehat{W}(b^{(k)}(\lambda))=e^{(k)}$ for each fixed
$k\in\mathbb{N}$, we obtain that
$$
\widehat{W}_n(z)=\sum_k\widehat{W}_n(b^{(k)}(\lambda))
=\sum_k\delta_{nk}=1\quad \textrm{for each} \quad n\in \mathbb{N}
$$
which shows that $\widehat{W}_n(z)=e\in c$ and hence $z\in
c^\lambda(\widehat{B})$. On the other hand, since
$c^\lambda(\widehat{B})$ and $c$ are the $BK-$spaces , Lemma
\ref{e5..1} implies the continuity of the matrix mapping
$A:c^\lambda(\widehat{B})\rightarrow c$. Thus, we have for every
$n\in\mathbb{N}$ that
$$
A_n(z)=\sum_kA_n(b^{(k)}(\lambda))=\sum_k\widehat{g}_{nk}.
$$
This show the necessity of (\ref{e5...2}).

Now, it follows by (\ref{e5..5}), (\ref{e5...1}) and (\ref{e5...2})
with Lemma\ref{e4..2} that
$\widehat{A}=(\widehat{g}_{n\widehat{k}})\in (c:c)$. So by
(\ref{e5.10}), (\ref{e5.11}) and relation (\ref{e5.15}) holds for
all $x\in c^\lambda(\widehat{B})$ and $y\in c$.

Finally, the necessity of (\ref{e5...0}) is immediately by
(\ref{e5.15}) since $Ax\in c$ and $\widehat{A}x\in c$. This
completes the proof.
\end{pf}
\begin{thm}
In order that $A=(a_{nk})\in (c^\lambda(\widehat{B}):c_0)$, it is
necessary and sufficient that (\ref{e5.10}), (\ref{e5.11}),
(\ref{e5..5})  hold and
\begin{eqnarray*}
&&\lim_na_n=0\\
&&\lim_n\widehat{g}_{nk}=0\quad\quad
 \textrm{ for each fixed}\ k\in \mathbb{N}\\
&&\lim_n\sum_k\widehat{g}_{nk}=0
\end{eqnarray*}
\end{thm}

\begin{pf}
This is obtained in the similar way used in the proof of
Theorem\ref{e5...5} with Lemma \ref{lm11} instead of Lemma
\ref{e4..2} and so we omit the detail.
\end{pf}
\begin{thm}\label{th16}
In order that $A=(a_{nk})\in (c_0^\lambda(\widehat{B}):c)$, it is
necessary and sufficient that  (\ref{e5.8}),
(\ref{e5..5}), (\ref{e5..6}) and (\ref{e5...1}) hold.
\end{thm}

\begin{pf}
This is an immediate consequence of Lemma \ref{e4.10}, Theorem
\ref{e4.004} and Theorem \ref{th13}(ii).
\end{pf}
\begin{thm}
In order that $A=(a_{nk})\in (c_0^\lambda(\widehat{B}):c_0)$, it is
necessary and sufficient that  (\ref{e5.8}), (\ref{e5..5}) and
(\ref{e5...1}) hold.
\end{thm}
\begin{pf}
This is an immediate consequence of Lemma \ref{lm12}, Theorem
\ref{e4.004} and Theorem \ref{th16}.
\end{pf}
\begin{thm}
Let $1<p<\infty$. $A=(a_{nk})\in
(\ell_p^\lambda(\widehat{B}):\ell_\infty)$ if and only if
 (\ref{e5.9})  and (\ref{e5..6}) hold and
\begin{eqnarray}\label{e5-3}
\sup_{n\in\mathbb{N}}\sum_{k=0}^{\infty}|\widehat{g}_{nk}|^q<\infty\qquad
(n\in \mathbb{N}).
\end{eqnarray}
\end{thm}
\begin{pf}
Suppose that $A=(a_{nk})$ satisfies the conditions  (\ref{e5.9}),
(\ref{e5..6}) and (\ref{e5-3}) and take any $x=(x_k)\in
\ell_p^\lambda(\widehat{B})$. Then, it is obvious that
$\{a_{nk}\}_{k\in \mathbb{N}}\in
\{\ell_p^\lambda(\widehat{B})\}^\beta$ for each $n\in\mathbb{N}$.
Hence, $A_{n}(x)=\sum_{k=0}^{\infty}a_{nk}x_k$ is convergent for all
$x\in \ell_p^\lambda(\widehat{B})$,i.e., $Ax$ exists. Also, it is
obvious by combining  Lemma\ref{e4..3} with the condition
(\ref{e5-3}) that $\widehat{A}=(\widehat{g}_{nk})\in
(\ell_p:\ell_\infty)$.

Now, we shall show that $Ax\in \ell_\infty$ for all $x\in
\ell_p^\lambda(\widehat{B})$. By using the relation (\ref{e5.13}) we
have that
\begin{eqnarray}\label{e5-4}
\sum_{k=0}^{\infty}a_{nk}x_k=\sum_{k=0}^{\infty}\widehat{g}_{nk}y_k.
\end{eqnarray}
By applying H\"{o}lder's inequality to (\ref{e5-4})we obtain that
\begin{eqnarray}\label{e5-5}
|(Ax)_n|&=&\bigg|\sum_{k=0}^{\infty}a_{nk}x_k\bigg|\leq\sum_{k=0}^{\infty}|
\widehat{g}_{nk}||y_k| \nonumber \\
&\leq& \|y\|_{\ell_p}\bigg(\sum_{k=0}^{\infty}|
\widehat{g}_{nk}|^q\bigg)^{1/q}\nonumber
\end{eqnarray}
which gives us by taking supremum over $n\in \mathbb{N}$ in this
last inequality that
$$
\sup_{n\in\mathbb{N}}|(Ax)_n|\leq
\|y\|_{\ell_p}\sup_{n\in\mathbb{N}}\bigg(\sum_{k=0}^{\infty}|
\widehat{g}_{nk}|^q\bigg)^{1/q}<\infty
$$
which implies that  $Ax\in \ell_\infty$ ,i.e., $A\in
(\ell_p^\lambda(\widehat{B}):\ell_\infty)$.

Conversely, suppose that $A=(a_{nk})\in
(\ell_p^\lambda(\widehat{B}):\ell_\infty)$. Thus, $Ax$ exists and
bounded for all $x\in \ell_p^\lambda(\widehat{B})$. Also,
$\{a_{nk}\}_{k\in
\mathbb{N}}\in\{\ell_p^\lambda(\widehat{B})\}^\beta$ for all
$n\in\mathbb{N}$ which is implies the necessity of (\ref{e5.9}) and
(\ref{e5..6}). If we define the sequences $\widehat{g}_n$ such that
$\widehat{g}_n=\{\widehat{g}_{nk}\}_{k\in \mathbb{N}}$ then we have
that
$$
\|\widehat{g}_n\|_{\ell_{q}}<\infty
$$
So, bearing in mind (\ref{e5..6}) one can easily obtain the relation
(\ref{e5-4}) by using the relation (\ref{e5.13}). On the other hand,
the sequences $a_n=\{a_{nk}\}_{k\in \mathbb{N}}$ define the
continuous linear functionals on the space
$\ell_p^\lambda(\widehat{B})$ as follows
$$
f_n(x)=\sum_{k=0}^{\infty}a_{nk}x_k;\qquad (n\in\mathbb{N}).
$$
Since the spaces $\ell_p^\lambda(\widehat{B})$ and $\ell_p^\lambda$
are linear isomorphic; we have that
$$
\|f_n\|=\|\widehat{g}_n\|_{\ell_q}.
$$
Also, since $\{f_n(x)\}_{n\in\mathbb{N}}\in \ell_\infty$; it is
obvious that
$$
\sup\|f_n(x)\|<\infty.
$$
Therefore, by using the Banach-Steinhaus Theorem \cite[see,
pp:1-2]{aw} we obtain that
$$
\sup_{n\in\mathbb{N}}\|f_n\|=\sup_{n\in\mathbb{N}}\|\widehat{g}_n\|_{\ell_q}
=\sup_{n\in\mathbb{N}}\bigg(\sum_{k=0}^{\infty}|
\widehat{g}_{nk}|^q\bigg)^{1/q}<\infty
$$
which shows that the necessity of (\ref{e5-3}). This step completes
the proof.
\end{pf}
\begin{thm}\label{thm19}
Let $1<p<\infty$. $A=(a_{nk})\in (\ell_p^\lambda(\widehat{B}):c)$ if
and only if (\ref{e5.9}), (\ref{e5..6}),(\ref{e5...1}) and
(\ref{e5-3}) hold.
\end{thm}
\begin{pf}
Suppose that the conditions (\ref{e5.9}),
(\ref{e5..6}),(\ref{e5...1}) and (\ref{e5-3}) hold. Then,\\
$\{a_{nk}\}_{k\in \mathbb{N}}\in
\{\ell_p^\lambda(\widehat{B})\}^\beta$ for each $n\in\mathbb{N}$.
Hence,the series  $A_{n}x=\sum_{k=0}^{\infty}a_{nk}x_k$ is
convergent for each $n\in\mathbb{N}$ and for all $x\in
\ell_p^\lambda(\widehat{B})$ which is implies that $Ax$ exists. By
using (\ref{e5...1}) we obtain that

$$
|\widehat{g}_{nk}|^q \rightarrow |\alpha_k|^q;\qquad
(n\rightarrow\infty).
$$
Since
$$
\bigg(\sum_{k=0}^{m}|
\widehat{g}_{nk}|^q\bigg)^{1/q}\leq\sup_{n\in\mathbb{N}}\bigg(\sum_{k=0}^{\infty}|
\widehat{g}_{nk}|^q\bigg)^{1/q}=M
$$
for all $m>0$. In this situation we see by passing to the limit in
last inequality as $n\rightarrow\infty$ that

\begin{eqnarray}\label{e5...6}
\bigg(\sum_{k=0}^{m}|
\alpha_{k}|^q\bigg)^{1/q}\leq\sup_{n\in\mathbb{N}}\bigg(\sum_{k=0}^{\infty}|
\widehat{g}_{nk}|^q\bigg)^{1/q}.
\end{eqnarray}
Because of (\ref{e5...6}) holds for all integer $m>0$; we have that
$$
\bigg(\sum_{k=0}^{\infty}| \alpha_{k}|^q\bigg)^{1/q}<\infty.
$$

We remember that $x=(x_k)$ and $y=(y_k)$ are associated  sequences
by the relation (\ref{e2.6}) where $y=(y_k)\in \ell_p$ for
$x=(x_k)\in\ell_p^\lambda(\widehat{B})$ .

Let $\varepsilon$ be any positive number. Then, there exists a
number $N$ such that
$$
\bigg(\sum_{k=N}^{\infty}|y_{k}|^q\bigg)^{1/q}<\frac{\varepsilon}{4M}.
$$
On the other hand, there is an integer $N_1$ such that
$$
\bigg|\sum_{k=0}^N\{\widehat{g}_{nk}-\alpha_k\}y_k\bigg|\leq\frac{\varepsilon}{2}
$$
whenever $n\geq N_1$. Therefore, we obtain
\begin{eqnarray*}
\bigg|\sum_{k=0}^\infty\{\widehat{g}_{nk}-\alpha_k\}y_k\bigg|&\leq&\bigg|\sum_{k=0}^N\{\widehat{g}_{nk}-\alpha_k\}y_k\bigg|+
\bigg|\sum_{k=N+1}^\infty\{\widehat{g}_{nk}-\alpha_k\}y_k\bigg| \\
&\leq&\frac{\varepsilon}{2}+\bigg(\sum_{k=N+1}^\infty\bigg[|\widehat{g}_{nk}|+|\alpha_k|\bigg]^q\bigg)^{1/q}
\bigg(\sum_{k=N+1}^\infty|y_k|^p\bigg)^{1/p} \\
&\leq& \frac{\varepsilon}{2}+2M\frac{\varepsilon}{4M}\\
&=&\varepsilon,
\end{eqnarray*}
which shows that
$$
\lim_{n\rightarrow\infty}\sum_{k=0}^\infty\widehat{g}_{nk}y_k=\sum_{k=0}^\infty\alpha_{k}y_k.
$$
Since
$$
\lim_{n\rightarrow\infty}(Ax)_n=\lim_{n\rightarrow\infty}\sum_{k=0}^\infty\widehat{g}_{nk}y_k=\sum_{k=0}^\infty\alpha_{k}y_k
$$
we have that $Ax\in c$ which is implies that
$A\in(\ell_p^\lambda(\widehat{B}):c)$.

The necessity part can be proved  by using the similar way in the
proof of Tehorem\ref{e5...5}, so we omit the detail.
\end{pf}
\begin{thm}
Let $1<p<\infty$. In order that  $A=(a_{nk})\in
(\ell_p^\lambda(\widehat{B}):c_0)$, it is necessary and sufficient
that (\ref{e5.9}), (\ref{e5..6}), (\ref{e5-3}) hold and
$$
\lim_n\widehat{g}_{nk}=0\quad\quad
 \textrm{ for each fixed}\ k\in \mathbb{N}
$$
\end{thm}
\begin{pf}
This result, can be proved similarly as the  proof of Theorem
\ref{thm19}.
\end{pf}












\end{document}